\documentclass[12pt]{article}
\usepackage[T2A]{fontenc}
\usepackage[cp1251]{inputenc}
\usepackage[russian]{babel}
\usepackage{amssymb,amsmath,latexsym,amsthm}
\usepackage{indentfirst}
\usepackage{verbatim}
\theoremstyle{plain}

\newtheorem{theorem}{Теорема}
\newtheorem*{theoremHB}{Теорема HB}
\newtheorem*{corollaryHB}{Следствие HB}
\newtheorem*{lemma}{Лемма}

\theoremstyle{remark}
\newtheorem{remark}{Замечание}

\DeclareMathOperator{\Hol}{Hol}
\newcommand{\CC}{\mathbb C}
\DeclareMathOperator{\psh}{psh}
\DeclareMathOperator{\sbh}{sbh}
\DeclareMathOperator{\dist}{dist}
\renewcommand{\leq}{\leqslant }
\renewcommand{\geq}{\geqslant }

\title{\begin{flushleft} {\normalsize УДК 517.55 + 517.576 + 517.987.1}   \end{flushleft}
{\large \bf Голоморфные миноранты плюрисубгармонических функций\footnote{Поддержано Российским фондом фундаментальных исследований, проект №\,16-01-00024а.}}}
\author{Т.\,Ю.~Байгускаров, Б.\,Н.~Хабибуллин}
\date{}
\begin{document}
\maketitle                                                                                                                                                    

Пусть $D$ --- область в $n$-мерном комплексном пространстве $\mathbb C^n$ с евклидовой нормой $|\cdot|$, $n\geq 1$. Через $\psh(D)$ и $\sbh (D)$, а также $\Hol (D)$  обозначаем  классы плюрисубгармонических и субгармонических, включая функцию $\boldsymbol{-\infty}$, тождественно равную $-\infty$ на $D$, а также голоморфных   в $D$  функций  соответственно; $B(z,r):=\bigl\{z'\in \CC^n \colon |z'-z|<r\bigr\}$ --- открытый шар с центром $z\in \CC^n$ радиуса $r>0$. Наиболее глубокий  результат о существовании ненулевой функции $f\in \Hol (D)$, в некотором смысле  минорирующей функцию    $\varphi \in \psh (D)$, ---  это доказанная Э.~Бомбьери [1, теорема существования]  методом Л.~Хермандера
\begin{theoremHB} [{[2; теорема 4.2.7]}] Пусть $D\subset \CC^n$ --- псевдовыпуклая область, $\varphi \in \psh (D)$, $z_0\in D$,
\begin{equation*}
I _{\varphi}^2(z_0,r) :=\int_{B(z_0,r)\cap D} e^{-2\varphi} \, d\lambda, \quad \text{где $\lambda$ --- мера Лебега на $\CC^n$}.
\end{equation*}
Тогда для каждого числа $a>0$ найдется функция $f_a\in \Hol (D)$, такая, что $f_a(z_0)=1$ и 
\begin{equation*}
	I_\varphi^2(z_0; f_a):= \int_D   \frac{\bigl| f_a(z)\bigr|^2e^{-2\varphi (z)}}{\bigl(1+|z-z_0|^2\bigr)^{n+a}}\, d \lambda (z)
	\leq \left(2+\frac{(1+r^2)^2(n+1)^2}{ar^{2n+2}} \right)	 \cdot I_{\varphi}^2(z_0,r).
\end{equation*}
\end{theoremHB}

Поскольку при $\boldsymbol{-\infty}\neq \varphi \in \psh (D)$ для  любой точки $z_0\notin \{z\in \CC^n \colon \varphi (z)=-\infty \}$ 
найдется число  $r>0$, для которого   $ I_{\varphi}^2(z_0,r) <+\infty$ [3, следствие 5.11], справедливо 
\begin{corollaryHB}  В условиях теоремы HB для каждого\/ $a>0$ найдутся\/ $z_0\in D$ и  ненулевая\/ $f_a \in \Hol (D)$, такие, что\/  $I_{\varphi}^2(z_0; f_a){<}+\infty$. В частности, для достаточно малого числа\/ $b>0$
\begin{equation}\label{fsup} 
	\log \bigl|bf_a(z)\bigr|\leq \sup_{z'\in B(z,r)} \varphi (z') + n\log \frac{1}{r}+(n+a)\log\bigl(1+|z|+r\bigr) \quad \text{для всех $z\in D$}
\end{equation}
 и  при всех\/ $0<r<\dist (z, \CC^n \setminus D)$, где $\dist$ --- функция евклидова расстояния.
 \end{corollaryHB}
Вывод  неравенства \eqref{fsup} из конечности $I_{\varphi}^2(z_0; f_a)$ элементарен, поскольку использует лишь субгар\-м\-о\-н\-и\-ч\-н\-ость функции  $|f_a|^2$ (см., например, [4, лемма 9.1]), но это неравенство находит весьма широкие применения в теории функций комплексных переменных и ее приложениях.  В случае $n=1$ и/или специальных областей $D\subset  \CC$ известны и более тонкие выводы, усиливающие  \eqref{fsup} и отчасти  учитывающие   как субгармоничность функций $\log |f_a|$ и $\varphi$, так и то, что $f_a\in \Hol (D)$:
\begin{description}
	\item[{\rm [1991 г.]}] {\it при\/ $D=\CC$ для\/ $r=1$ в правой части неравенства\/ \eqref{fsup} можно 
	оставить лишь\/} $\sup_{z'\in B(z,1)} \varphi (z')$    [5, доказательство основной теоремы, 2а)--2б)], [4, предложение 9.2];
	\item[{\rm [1992 г.]}] {\it 	 для  непрерывной ограниченной сверху функции\/ $r\colon D\to (0,+\infty)$ в области $D\subset \CC$,  когда $r(z)<\dist (z, \mathbb C\setminus D)$, $\sup \bigl\{r(z)/r(z')\colon z,z'\in D,  |z-z'|\leq r(z)+r(z')\bigr\}<+\infty$ и $-\log r \in \sbh(D)$,   оценки\/ \eqref{fsup} можно усилить усреднением функции\/ $\varphi$, заменив правую часть на
\begin{equation}\label{EO}
	\frac{1}{2\pi} \int_{-\pi}^{\pi} \varphi \bigl(z+r(z)e^{i\theta}\bigr) \,d \theta+\log \frac{1}{r(z)} +4\log \bigl(1+|z|\bigr) \qquad \text{{\rm {\large(}см.  [6, лемма]{\large)}}};
\end{equation}
}
\item[{\rm [2015 г.]}] {\it при\/ $D=\CC$  для любого числа $r\in (0,1]$ при достаточно малом $b>0$ правую часть в\/ \eqref{fsup}  можно  заменить на усреднение 
$\frac{1}{2\pi} \int\limits_{-\pi}^{\pi} \varphi (z+re^{i\theta}) \, d \theta$ {\rm [7, теорема~2]}.}
\end{description}
 \begin{theorem}\label{th1} Для псевдовыпуклой области $D\subset \CC^n$ при $\boldsymbol{-\infty} \neq \varphi \in \psh (D)$ для любого числа $a>0$ 
  найдется ненулевая функция $f\in \Hol (D)$, удовлетворяющая  оценке  
\begin{equation}\label{est:main}
	\log \bigl|f(z)\bigr|\leq \frac{n!}{\pi^nr^{2n}} \int_{B(z,r)} \varphi  \, d \lambda + n\log \frac{1}{r}+(n+a)\log\bigl(1+|z|+r\bigr) 
\end{equation}
 для всех $z\in D$ и всех\/ $r$, $0<r<\dist (z, \CC^n \setminus D)$, где\/ ${n!}/{\pi^nr^{2n}}={1}/{\lambda(B(0,r))}$, т.\,е. первое слагаемое  в правой части неравенства\/ \eqref{est:main} --- усреднение\/ $\varphi$ в шаре\/ $B(z,r)$ по мере Лебега\/ $\lambda$.
\end{theorem}
\begin{proof} По первому утверждению cледствия HB найдутся точка $z_0\in D$,  ненулевая функция $f_a \in \Hol (D)$ и число $C_a>0$, для которых  
 $I_{\varphi}^2(z_0; f_a){\leq }C_a$, т.\,е.
\begin{equation*}
	\int_D   \bigl| f_a\bigr|^2e^{-2\psi}\, d \lambda  \leq C_a, \quad 
\psi(z):= \varphi (z)+\frac{1}{2}(n+a)\log \bigl(1+|z-z_0|^2\bigr), \quad z\in D.	
\end{equation*}
\begin{lemma} Пусть $u \in \sbh \bigl(B(z,r')\bigr)$ и  $r'>r>0$. Тогда для любой функции $v$, интегрируемой по мере Лебега $\lambda$ на $B(z,r)$,  справедливы неравенства
\begin{multline*}
	\exp \Bigl(u(z)-\frac{1}{\lambda\bigl(B(0,r)\bigr)} \int_{B(z,r)} v \, d \lambda \Bigr)\leq  	\exp \Bigl(\frac{1}{\lambda\bigl(B(0,r)\bigr)} \int_{B(z,r)}(u- v)\, d \lambda \Bigr)
	\\ 
	\leq  	\frac{1}{\lambda\bigl(B(0,r)\bigr)} \int_{B(z,r)} e^{u-v} \, d \lambda.
\end{multline*}
\end{lemma}  
Первое неравенство здесь следует из условия $u\in \sbh \bigl(B(z,r')\bigr)$, а второе ---  из интегрального неравенства Йенсена для выпуклой функции $\exp$ [8, теорема 2.1].

Применим лемму к  $u{:=}2\log |f_a|\in \psh (D)\subset \sbh (D)$ и $v{:=}2\psi$ для  произвольных точки $z\in D$ и  $r$, $0<r<\dist \bigl(z, \CC^n\setminus D\bigr)$ (ниже $V_{2n}$ --- объем единичного шара в $\CC^n$):
 \begin{multline*}
2\log \bigl|f_a(z)\bigr| {\leq} 2\,\frac{1}{\lambda\bigl(B(0,r)\bigr)} \int_{B(z,r)} \psi \, d \lambda+ \log \frac{1}{\lambda\bigl(B(0,r)\bigr)}+\log  \int_{B(z,r)} e^{2\log |f_a|-2\psi} \, d \lambda\\ {\leq} 2\,\frac{1}{\lambda\bigl(B(0,r)\bigr)} \int_{B(z,r)} \varphi \, d \lambda+
\frac{n+a}{\lambda\bigl(B(0,r)\bigr)} \int_{B(z,r)} \log \bigl(1+|z'-z_0|^2\bigr) \, d \lambda (z')+2n\log \frac{1}{r} -\log V_{2n}\\+\log  \int_{B(z,r)} |f_a|e^{-2\psi} \, d \lambda {\leq} 2\,\frac{1}{\lambda\bigl(B(0,r)\bigr)} \int_{B(z,r)} \varphi \, d \lambda+(n+a) \log \bigl(1+\bigl(|z-z_0|+r\bigr)^2\bigr)  +2n\log \frac{1}{r}\\-\log V_{2n}+\log C_a \leq 2\,\frac{1}{\lambda\bigl(B(0,r)\bigr)} \int_{B(z,r)} \varphi \, d \lambda+2n\log \frac{1}{r} +
2(n+a) \log \bigl(1+|z|+r\bigr) \\+2(n+a)\log \bigl(1+|z_0|\bigr) -\log V_{2n}+\log C_a.
 \end{multline*}
Деление неравенств на $2$  при $f:=bf_a$ с числом $b:=(1+|z_0|)^{-n-a}\sqrt{V_{2n}/ C_a}$ дает \eqref{est:main}.  
\end{proof}
\begin{remark} При $n=1$ правая часть в \eqref{est:main} всегда, причем зачастую и существенно, меньше и лучше, чем \eqref{EO}. Этого нельзя сказать о двух других предшествующих результатах от [1991~г.] и [2015~г.]. Но нам удалось установить, что при $\boldsymbol{-\infty} \neq \varphi \in \sbh (\CC)$  для любого сколь угодно большого числа $N\geq 0$  найдется ненулевая функция $f\in \Hol (\CC)$, для которой
\begin{equation*}
 \log \bigl|f(z)\bigr| \leq \frac{1}{2\pi} \int_0^{2\pi} \varphi\left(z+\frac{1}{\bigl(1+|z|\bigr)^N}\, e^{i\theta}\right) \, d \theta 
\end{equation*}
для всех $z\in \CC$. Доказательство  этого [9, следствие 2] и даже более тонких результатов для голоморфных функций одной комплексной переменной см. в [9].
\end{remark}

Пусть  $\mathcal M_{\rm c}^+(D)$ --- класс всех положительных борелевских мер с компактным носителем в области $D\subset \CC^n$, а $\mathcal M_{\rm ac}^+(D)$ --- подкласс всех абсолютно непрерывных относительно  $\lambda$ мер в нем. Мера  $\mu\in \mathcal M_{\rm c}^+(D)$ --- называется {\it выметанием меры\/} $\mu_0\in \mathcal M_{\rm c}^+(D)$ относительно $\psh (D)$, если 
\begin{equation*}
	\int \varphi \, d \mu_0\leq  \int \varphi \, d \mu \quad \text{{\it для всех\/} $\varphi \in \psh (D)$\quad  (пишем $\mu_0\prec_{\psh(D)}\mu$).}
\end{equation*}
\begin{theorem}\label{th2} Пусть область $D\subset \CC^n$ псевдовыпукла,  $\mu_0\in \mathcal M_{\rm c}^+(D)$, а   $F\colon D\to [-\infty, +\infty]$ --- интегрируемая по  $\mu_0$ и локально интегрируемая в $D$ по мере Лебега $\lambda$ функция. Если 
\begin{equation}\label{infF}
-\infty <\inf \Bigl\{ \,\int F\, d \mu\colon \mu_0\prec_{\psh(D)} \mu \in \mathcal M_{\rm ac}^+(D) \Bigr\}	,
\end{equation}
то для любых числа $a>0$ и непрерывной функции $d\colon D\to (0,+\infty)$, удовлетворяющей условию   $d(z)<\dist (z, \CC^n\setminus D)$ при всех $z\in D$, для  <<переменного>> усреднения  
\begin{equation}\label{dfFd}
F^{*d}(z):=\frac{1}{\lambda(B(0,d(z)))} \int_{B(z,d(z))}F\, d\lambda, \quad z\in D,
\end{equation} 
найдется ненулевая функция $f\in \Hol (D)$, для которой 
\begin{equation}\label{fF}
	\log \bigl|f(z)\bigr|\leq \frac{1}{\lambda\bigl(B(0,r)\bigr)} \int_{B(z,r)}F^{*d}\, d\lambda + n\log \frac{1}{r}+(n+a)\log\bigl(1+|z|+r\bigr) 
\end{equation}
для всех\/ $z\in D$ и\/ $0<r<\dist (z, \CC^n \setminus D)$.  
\end{theorem}
\begin{proof} Условие \eqref{infF} означает, что для функции $F^{*d}$ из \eqref{dfFd} найдется миноранта $\varphi\leq F^{*d}$ на $D$, $\boldsymbol{-\infty}\neq \varphi\in \psh(D)$, если воспользоваться теоремой 7.1 из  [4] в пространстве $\CC^n$, отождествленном с $2n$-мерным вещественным пространством, c $H:=\psh (D)\subset \sbh(D)$ и $F^{(d)}:=F^{*d}$. Теперь достаточно применить  теорему \ref{th1} к  миноранте $\varphi\in \psh (D)$. 
\end{proof}
\begin{remark} Условие \eqref{infF}, достаточное для существования ненулевой голоморфной миноранты для переменного усреднения функции $F$ в форме \eqref{fF} очень близк\'о к необходимому [4, предложение 7.1], если  отвлечься от <<логарифмических добавок>>  и усреднений функции $F$ в правой части неравенства \eqref{fF}.
\end{remark}
\begin{center}
{\sc Литература}
\end{center}

\noindent
[1] E.~Bombieri, Invent. Math., {\bf 10} (1970), 248--263. 

\noindent
[2] L.~H\"ormander, {\it Notions of Convexity,\/} Birkh\"aser, Boston, 1994. 

\noindent
[3]  П.~Лелон, Л.~Груман, {\it Целые функции многих комплексных переменных,\/} Мир, М., 1989. 

\noindent
[4] Б.\,Н.~Хабибуллин,  Изв. РАН, cер. матем., {\bf 65}:5   (2001), 167--190. 

\noindent
[5]  Б.\,Н.~Хабибуллин,  Изв. АН СССР, cер. матем. {\bf 55}:5   (1991). 1101--1123. 

\noindent
[6]  О.\,В.~Епифанов, Матем. заметки, {\bf 51}:1 (1992). 83--92. 

\noindent
[7]  Т.\,Ю.~Байгускаров, Г.\,Р.~Талипова, Б.\,Н.~Хабибуллин, Алгебра и анализ, 2016 принята к печати.

\noindent
[8] У.~Хейман, П.~Кеннеди, {\it Субгармонические функции,\/} Мир, М., 1980.

 \noindent
[9]  Б.\,Н.~Хабибуллин, Т.\,Ю.~Байгускаров, Матем. заметки, 2016 принята к печати.
\end{document}